\documentclass[lettersize,journal]{article}

\usepackage{amsmath,amsfonts}
\usepackage{algorithmic}
\usepackage{algorithm}
\usepackage{array}
\usepackage[caption=false,font=normalsize,labelfont=sf,textfont=sf]{subfig}
\usepackage{textcomp}
\usepackage{stfloats}
\usepackage[margin=1in]{geometry}
\usepackage{url}
\usepackage{verbatim}
\usepackage{graphicx}
\usepackage{cite}

\usepackage{hyperref}
\hypersetup{
	colorlinks=true,
	breaklinks=true,
	urlcolor=blue,
	linkcolor=blue,
	bookmarksopen=true,
	citecolor=blue,
}

\usepackage{amsthm}
\theoremstyle{definition}

\theoremstyle{remark}

\newtheorem{lemma}{Lemma}[section]

\newtheorem{theorem}{Theorem}[section]

\usepackage{authblk}

\title{A limiting model for a low {R}eynolds number swimmer with {$N$} passive elastic arms}

\author[1]{François Alouges}
\author[2]{Aline Lefebvre-Lepot}
\author[2]{Jessie Levillain}

\affil[1]{Centre Borelli, ENS Paris-Saclay, CNRS, Université Paris-Saclay, 91190 Gif-sur-Yvette, France}
\affil[2]{ CMAP, CNRS, École polytechnique, Institut Polytechnique de Paris, 91120 Palaiseau, France}

\date{}

\begin{document}

\maketitle

\begin{abstract}
We consider a low Reynolds number artificial swimmer that consists of an active arm followed by $N$ passive springs separated by spheres. This setup generalizes an approach proposed in Montino  and  DeSimone, Eur.  Phys.  J. E, vol. 38, 2015. We further study the limit as the number of springs tends to infinity and the parameters are scaled conveniently, and provide a rigorous proof of the convergence of the discrete model to the continuous one. Several numerical experiments show the performances of the displacement in terms of the frequency or the amplitude of the oscillation of the active arm. \end{abstract}

\section{Introduction}

As stated by Purcell's {\emph{Scallop Theorem}} \cite{purcell_life_1977}, reciprocal shape changes in a swimmer never leads to a net displacement of the system in a low Reynolds number setting. Indeed, a microscopic scallop opening and closing its valve would be completely unable to swim, due to negligible inertial forces in this situation \cite{childress_1981}. Several simple mechanisms have then been introduced (see e.g. \cite{lauga_hydrodynamics_2009}) to overcome this obstruction, most of them using two degrees of freedom in order to create closed curves with nonzero surface in the shape space of the swimmer.\\

One of the simplest mechanisms introduced in the literature is probably Najafi and Golestanian's three-sphere swimmer \cite{najafi_simple_2004}, which consists in three spheres linked by two extensible arms of negligible thickness, moving along a single direction. This model is much simpler than Purcell's original three-link swimmer \cite{purcell_life_1977}, or Purcell's rotator \cite{dreyfus_purcells_2005}, as there is no rotational motion involved. This swimmer has two degrees of freedom, activated periodically in time with a phase lag in order to produce the loop. Both Purcell's and Najafi and Golestanian's swimmers have been extensively studied in \cite{alouges_optimally_2013, desimone_computing_2012, alouges_num, alouges_optimal_2008, alouges_energy-optimal_2019, alouges_optimal_2009}.\\
As an extension of this three-sphere swimmer, Montino and DeSimone then introduced a three-sphere swimmer with a passive elastic arm \cite{montino_three-sphere_2015}. This adaptation has only one degree of freedom, which is the length of the non-elastic arm. Thanks to a resonant effect at natural frequency of the system (which depends on the viscosity of the fluid, the masses and the spring constant), an out-of-phase oscillation of the spring is created, which ultimately leads to a net motion of the swimmer. However, at very low or very high frequency, no net motion is possible after a stroke.
Having this passive elastic arm also confines net motion to only one direction on the swimming axis, swimming direction is thus limited, and the swimmer can only move with its passive arm ahead. This was also denoted by Passov \cite{passov_dynamics_2012}, when looking at Purcell's three-link swimmer with a passive elastic tail.\\

In this paper, Montino and DeSimone's swimmer is extended by adding a large number $N$ of passive elastic arms to their one-dimensional swimmer, thus turning it into an $(N+2)$-sphere swimmer. This simple swimmer then leads to a limit model with an elastic tail resembling a one-dimensional flagella along which compressive waves propagate.\\

The paper is organized as follows. In sect. \ref{sec:Narms}, we describe the $N$-spring swimmer, and its equations of motion, before looking at the limit model, when the number of springs tends to infinity, in sect. \ref{sec:cont}. We prove the convergence of the discrete model to the continuous one in sect. \ref{sec: dtoc}, using the fact that it is found to be a non-conventional mass lumping discretization of the limit model. Sect. \ref{sec:dx} introduces two formulas in order to compute the net displacement of both swimmers, discrete and continuous.  Finally, in sect. \ref{sec:num} we study numerically the movement and displacement of our swimmer depending on various system parameters, in order to find optimal swimming parameters to obtain the largest net displacement possible.

\section{Problem's formulation and study: $N$-spring discrete model and its continuous limit}\label{sec:Narms}

\begin{figure}[!h]
    \centering
    \includegraphics[width = 0.75\textwidth]{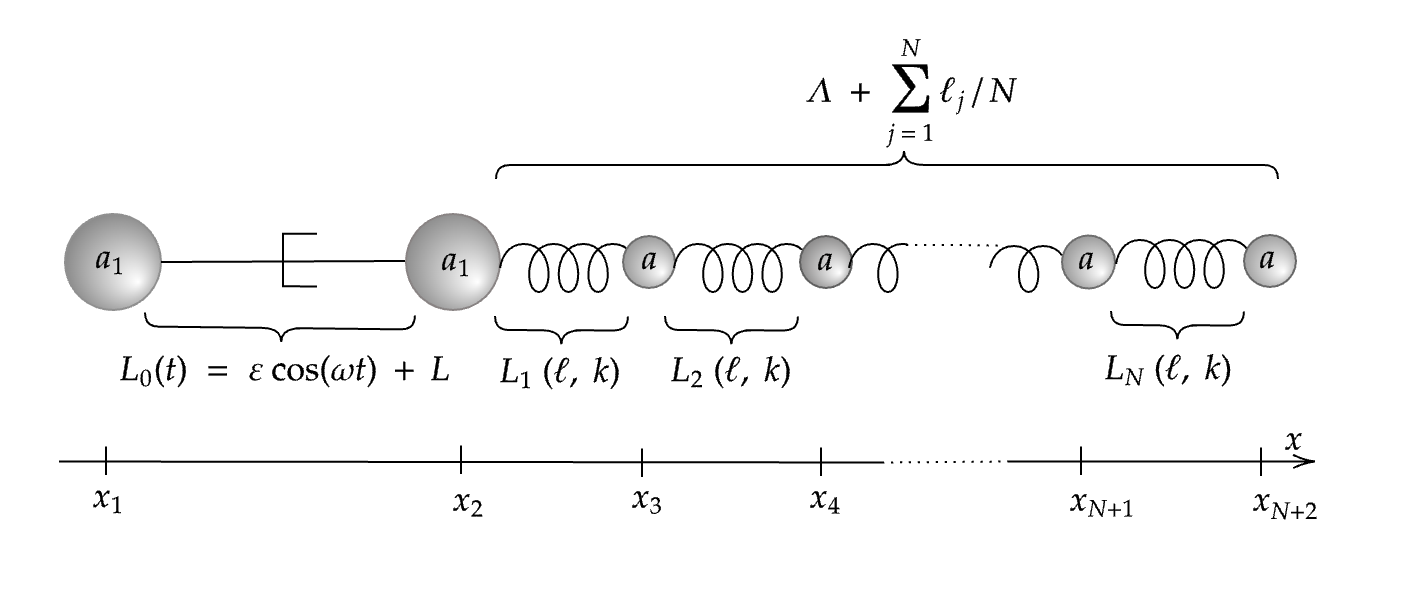}
    \caption{Low Reynolds number swimmer with $N$ elastic arms.}
    \label{fig:Nres}
\end{figure}

The swimmer studied in this paper is an extension of the three-sphere swimmer with a single passive elastic arm \cite{montino_three-sphere_2015}, to a swimmer with $N+2$ spheres and $N$ passive elastic arms, presented in figure \ref{fig:Nres}. The first arm of this artificial swimmer is a rod of negligible thickness, surrounded by two spheres of radius $a_1$. This arm has a prescribed periodic movement around a length at rest $L$, of the form $L_0(t) = L(1+\tilde \varepsilon \cos(\omega t))$ where $\tilde \varepsilon\in [0,1)$ is a non-dimensional parameter. $\tilde \varepsilon<1$ so that the active arm always has a positive length. We define $\varepsilon$ as  $\varepsilon = L \tilde \varepsilon$. The rest of the swimmer has a total length  at rest $\Lambda$ that does not depend on $N$. In order to keep a constant length and have an elastic force that does not depend on $N$, all the other spheres have a radius $a = \tilde a/N$, the springs have each a rest length $h = \Lambda/N \gg a$, and an elastic constant $k = \tilde kN$, with $\tilde k$ and $\tilde a$ prescribed and independent of $N$. 

If the swimmer is able to control the length of the front rod with the prescribed periodic function $L_0(t)$, the length of the $N$ remaining springs are governed by the balance of viscous and elastic forces.
At any time $t$, the length $L_j(t)$ of the $j$-th arm, $j \geq 1$ is written as $L_j(t) = \frac{\ell_j(t)}{N}+h$. 
Let us then denote by $\mu$ the fluid viscosity, $f_j^F$ and $f_j^R$ the hydrodynamic and elastic forces on the $j$-th sphere. We also call $x_j$ the coordinate of its center, so that $V_j = \dot x_j$ is the velocity of the $j$-th sphere. The geometry of the system entails $\dot L_j = V_{j+2}-V_{j+1}$ for all $j=0,\, \dots,\, N $.\\

In order to effectively swim, our $N$-spring swimmer undergoes periodic harmonic but non-reversible deformations, just like the original swimmers from Najafi and Golestanian \cite{najafi_simple_2004}, and Montino and DeSimone \cite{montino_three-sphere_2015}. However, due to the geometry, we expect a wave to propagate along the tail. This is the behaviour of this wave that we aim at describing in the remainder of the paper.

\subsection{First approximations}

In a first approximation, we consider the case where the hydrodynamic force on the $j$-th sphere only depends on the speed of that same sphere, and neglect interactions between spheres. This leads to the following set of equations on (fluid) forces and velocities:

\begin{equation}
\left\{\begin{array}{ll}
        f_j^F = -6\pi \mu a V_j \text{ for } j \geq 3,\vspace{6pt} \\
        f_j^F = -6\pi \mu a_1 V_j \text{ for } j =1,2.
    \end{array}
    \right.
    \label{eq:forcesfl}
\end{equation}

The elastic forces on each sphere can be written as:

\begin{equation}
\left\{\begin{array}{ll}
f_{2}^R &\displaystyle = k(L_{1}-h) = k\frac{\ell_{1}}{N}\vspace{6pt}\\
        f_j^R &= k\big((L_{j-1}-h) - (L_{j-2}-h)\big),\vspace{6pt}\\
        &\displaystyle  = k\frac{\ell_{j-1}-\ell_{j-2}}{N}  \quad \text{ for } 3 \leq j \leq N+1\vspace{6pt}\\
        f_{N+2}^R &\displaystyle = -k(L_{N}-h) = -k\frac{\ell_{N}}{N}.
    \end{array}
    \right.    
    \label{eq:forcesel}
\end{equation}

At low Reynolds number, inertial forces are negligible. This, together with the fact that the artificial swimmer is self-propelled, gives:

\begin{equation}
\left\{\begin{array}{ll}
        f_1^F + \dots + f_{N+2}^F = 0,\vspace{6pt} \\
        f_j^R + f_j^F =0 \text{ for } j \geq 3.
    \end{array}
    \right.
    \label{eq:pfd}
\end{equation}

Using \eqref{eq:forcesfl}, \eqref{eq:forcesel} and \eqref{eq:pfd}, we obtain the expression of fluid forces on each sphere with respect to the length of the adjacent arms. In particular, for the first two spheres:
\begin{equation}
    \left\{\begin{array}{ll}
        f_1^F - f_2^F = 6 \pi \mu a_1 (V_2 - V_1) = 6 \pi \mu a_1 \dot L_0,\vspace{6pt}\\
        f_1^F + f_2^F  = f_3^R + \dots + f_{N+2}^R =  -k\ell_1/N,
    \end{array}
\right.
\end{equation}
which finally leads to:
\begin{equation}
\left\{\begin{array}{ll}
        f_1^F = \frac{1}{2}(+6\pi \mu a_1 \dot L_0 - \tilde k \ell_1),\vspace{6pt}\\
        f_2^F =  \frac{1}{2}(-6\pi \mu a_1 \dot L_0 - \tilde k \ell_1).
    \end{array}
\right.     
\label{eq:ff}
\end{equation}

\subsection{Movement of the spheres}\label{sec:mvt}

In order to write the equations governing the system, we use equations (\ref{eq:forcesfl}-\ref{eq:ff}) to find ODEs on the elongation $l_j(t)$ of the $j$-th arm, for $j \geq 1$. We first consider the case $j \geq 2$. Writing $\dot L_j = V_{j+2}-V_{j+1} = \frac{1}{6\pi\mu a}(f^R_{j+2}-f^R_{j+1})$, one deduces
\begin{equation}
    \displaystyle \dot \ell_j = \Lambda^2 K \frac{\ell_{j-1} - 2 \ell_j + \ell_{j+1}}{h^2}, \, 2\leq j \leq N,\\
    \label{eq:elli}
\end{equation}
where we have added a fictitious variable 
\begin{equation}
    \ell_{N+1}=0\,,
    \label{eq: lN+2}
\end{equation}
and with $K = \displaystyle \frac{\tilde k}{6 \pi \mu \tilde a}$ . 

To determine the equation for the first elastic arm,  we use the fact that $\displaystyle \dot L_1 = V_3 - V_2= -\frac{1}{6\pi\mu a}f^F_{3}+\frac{1}{6\pi\mu a_1}f^F_{2}$ to obtain, using equations \eqref{eq:forcesel} and \eqref{eq:ff}: 

\begin{equation}
   h \dot \ell_1 = \Lambda^2 K\frac{\ell_2 - \ell_1}{h}-\frac{ \Lambda K \tilde a}{2 a_1}\ell_1 -\frac{\Lambda}{2} \dot L_0.
    \label{eq: l2}
\end{equation}

We can easily verify that the ODE problem (\ref{eq:elli},\ref{eq: lN+2},\ref{eq: l2}) is well-posed using Cauchy-Lipschitz theorem, and provides a unique solution $(\ell_j(t))_{1\leq j \leq N+1}$ for any initial configuration.

Seeking for periodic (complex) solutions to equation \eqref{eq:elli} leads to 
\begin{equation}
\ell_j(t) = (\alpha_d \gamma_+^{j-1} + \beta_d \gamma_-^{j-1})e^{i \omega t},
\label{eq:Disc_solexpl1}
\end{equation}
where $\alpha_d, \beta_d \in \mathbb{C}$ and
\begin{equation}
\gamma_\pm = \frac{i/(K_{\omega} N^2) +2 \pm \sqrt{\Delta}}{2} 
\label{eq:gamma}
\end{equation}
and $\displaystyle \Delta = \frac{-1}{K_{\omega}^2N^4} + \frac{4i }{K_{\omega}N^2}$, where $\displaystyle K_{\omega} = \frac{K}{\omega} = \frac{\tilde k}{6 \pi \mu \tilde a \omega}$ is an adimensional number. Notice that $\vert \gamma_+ \vert >1$ while $\vert \gamma_- \vert <1$. The constants $\alpha_d$ and $\beta_d$ may be determined through the boundary conditions. Namely assuming, from the linearity of the problem, $\ell_1=b_d e^{i\omega t}$, with $b_d\in \mathbb{C}$ and recalling $l_{N+1}=0$ enables us to write
\begin{equation}
\left\{\begin{array}{ll}
        \displaystyle\ell_1 (t) = b_d e^{i \omega t} = e^{i\omega t}(\alpha_d + \beta_d),\vspace{6pt}\\
        \displaystyle\ell_{N+1} (t) = e^{i\omega t}(\alpha_d \gamma_+^{N} + \beta_d \gamma_-^{N}) = 0\,,
    \end{array}
\right.
\label{eq:boundary}
\end{equation}
to finally obtain
 \begin{equation}
    \alpha_d = \frac{-\gamma_-^N b_d}{(\gamma_+^N - \gamma_-^N)},\quad 
    \beta_d = \frac{\gamma_+^N b_d}{(\gamma_+^N - \gamma_-^N)}.
    \label{eq:ab}
\end{equation}

Then, we use \eqref{eq: l2} to determine $b_d$:
\begin{equation}
b_d = - \frac{\varepsilon i /2}{\displaystyle i/N + N K_{\omega}(1-z_d) + K_{\omega}\frac{\tilde a}{2 a_1}},
\label{eq: bd}
\end{equation}
where $\displaystyle z_d = \frac{\gamma_+^N\gamma_- - \gamma_-^N\gamma_+}{\gamma_+^N - \gamma_-^N}$ .

\subsection{Limit model with an infinite number of springs}\label{sec:cont}

\begin{figure}[!ht]
    \centering
    \includegraphics[width = 0.75 \textwidth]{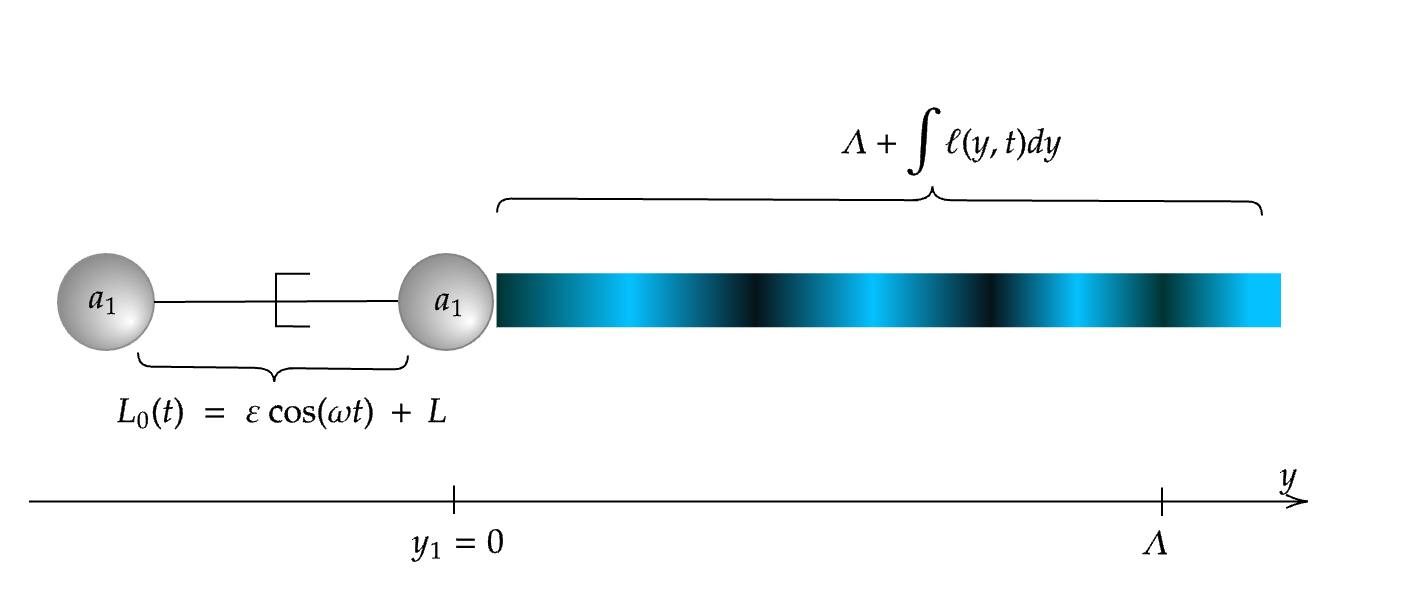}
    \caption{Continuous model of the low-Reynolds-number elastic swimmer. Color variations in the tail indicate compression and expansion of the swimmer.}
    \label{fig:resC}
\end{figure}

As we increase the number of springs in our swimmer, a limit model arises, with an elastic-like tail, as shown in figure \ref{fig:resC}. This elastic tail compresses and dilates itself in the same way that the springs do, following the active arm, in order to create a global displacement of our swimmer.  

Equations \eqref{eq:elli} and \eqref{eq: l2} can be viewed as a finite element discretization of a PDE, which describes the continuous version of our swimmer. Limit expressions for this PDE model are formally derived throughout this section while the convergence of the $N$-spring model to the continuous model will be proven in Sect. \ref{sec: dtoc}.

First, as $h \to 0$ ($N \to \infty$), $\displaystyle \frac{\ell_{j-1} -2 \ell_j +\ell_{j+1}}{h^2}$ formally converges to a second order derivative. More precisely, we introduce a new space variable $y_j = (j-1)h$ for $1\leq j \leq N+1$. The points $y_j$ are equally spaced and thus different from the previous $x_j$. Since $y_1 = 0$, the $y$ variable can be seen as a local space coordinate attached to the second sphere, and we assume $\ell(y_{j}) = \ell_{j}$ for a smooth enough function $\ell$. Passing to the formal limit in \eqref{eq:elli} leads to a heat equation:
\begin{equation} 
    \partial_t \ell(y,t) = K\Lambda^2 \partial_{yy} \ell(y,t), \quad \forall (y,t) \in [0, \Lambda] \times \mathbb{R}^{\star}_+.
    \label{chaleur}
\end{equation}

Concerning the boundary conditions, we first notice that $\ell_{N+1}=0$ leads to $\ell(\Lambda,t) = 0$ for all $t>0$.
As $h \to 0$, the equation \eqref{eq: l2} on $\ell_1$ formally becomes a Fourier-type boundary condition:
\begin{equation*}
    \Lambda^2 K \partial_y \ell(0,t) - \Lambda K \frac{\tilde a}{2a_1}\ell(0,t) = \frac{\Lambda}{2}\dot L_0(t), \quad \forall t>0.
\end{equation*}

Therefore, we finally obtain the following continuous problem: 

Find $\ell \in \mathcal{C}^2([0, \Lambda] \times \mathbb{R}^*_+)$ such that $\forall (y,t) \in (0, \Lambda) \times \mathbb{R}^{\star}_+,$
  \begin{equation}
\left\{\begin{array}{l}
         \displaystyle \partial_t \ell(y,t) - \Lambda^2 K\partial_{yy} \ell(y,t)=0, \vspace{6pt}\\
        \displaystyle \Lambda^2 K \partial_y \ell(0,t) - \Lambda K \frac{\tilde a}{2a_1}\ell(0,t) = \frac{\Lambda}{2}\dot L_0(t), \vspace{6pt}\\
        \displaystyle \ell(\Lambda,t) = 0.
    \end{array}
\right.       
\label{eq : continue}
  \end{equation}

\subsection{Well-posedness of the problem}

Equation \eqref{eq : continue} belongs to the class of problem for which the classical theory of parabolic equations applies. Namely, calling
$$
\mathcal{V} = \left\{u\in H^1((0,\Lambda)) | u(\Lambda)=0\right\}\,,
$$
which is a Hilbert space with the scalar product $(u,v)_{\mathcal{V}}=\int_0^\Lambda \partial_y u \,\partial_y v\,dy$, the variational formulation reads:

Let $T>0$, find $\ell(y,t)\in L_t^\infty(0,T;L_y^2((0,\Lambda)))\cap L_t^2(0,T;\mathcal{V})$ such that for all $t\in (0,T)$ and for all $v\in \mathcal{V}$
\begin{eqnarray}
&&\frac{d}{dt}\int_0^\Lambda \ell v\,dy + \Lambda ^2 K\int_0^\Lambda \partial_y \ell \,\partial_y v\,dy\label{eq:FV} \\
&&\hspace*{2cm}+ \frac{\Lambda K\tilde a}{2 a_1}\ell(0,t)v(0) = -\frac{\Lambda}{2}\dot L_0(t)v(0)\nonumber
\end{eqnarray}
with $\ell(y,0)=\ell_0(y)\in L^2((0,\Lambda))$ a given initial data.

Defining the bilinear form $\kappa$ in $\mathcal{V} \times \mathcal{V}$ as:
\begin{equation}
    \kappa : (u,v) \mapsto \Lambda^2 K \int_0^\Lambda \partial_y u(y)\partial_y v(y)\,dy + \frac{\Lambda K \tilde a}{2 a_1}u(0)v(0),
    \label{eq: Abilin}
\end{equation}
which is symmetric and coercive on $\mathcal{V}$, well-posedness of the problem (\ref{eq:FV}) follows from standard results on parabolic equations (see e.g. \cite{raviart_introduction}). Moreover, it is well known that the solution $\ell(\cdot,t)$ is of class $\mathcal{C}^\infty([0,\Lambda])$ for any time $t>0$.

\subsection{Analytical periodic solutions}

Let us now solve the system \eqref{eq : continue} using the following ansatz $\ell(y,t) = \underline{\ell}(y)e^{i \omega t}$.
From \eqref{chaleur} we deduce the following equation:
\begin{equation}
    i \underline{\ell}= \Lambda^2 K_{\omega} \partial_{yy} \underline{\ell}\,.
    \label{eq:lbar}
\end{equation}

The characteristic polynomial associated to \eqref{eq:lbar} has two roots, $\displaystyle r := \frac{1+i}{\Lambda \sqrt{2 K_{\omega}}}$ and $-r$, which leads to the following solutions:
\begin{equation}
    \underline{\ell}(y) = \alpha e^{ry}+\beta e^{-ry},
\label{eq:Cont_solexpl1}
\end{equation}
with $\alpha, \, \beta \in \mathbb{C}$.

We then determine $\alpha$ and $\beta$ using boundary conditions:
$$\left\{\begin{array}{ll}
        \displaystyle -(\alpha+\beta) \frac{\tilde a}{2 a_1} + \Lambda r(\alpha  - \beta) = \frac{i \varepsilon}{2 K_{\omega}},\vspace{6pt}\\
        \displaystyle \alpha e^{r\Lambda}+ \beta e^{-r\Lambda} = 0,
    \end{array}
\right.         
$$
i.e.,
\begin{equation}
    \left\{\begin{array}{ll}
    \displaystyle \alpha = \frac{i \varepsilon }{\displaystyle 2 K_{\omega} \big( \frac{\tilde a}{2a_1} (e^{2r\Lambda}-1) +\Lambda r(e^{2r\Lambda}+1)\big)},\vspace{6pt}\\
    \displaystyle \beta = -e^{2r\Lambda} \alpha.
        \end{array}
\right.  
\label{eq:Cont_solexpl2}
\end{equation}

We notice that $\displaystyle r \Lambda = \frac{1+i}{\sqrt{2 K_{\omega}}}$ only depends on $K_\omega$.

\section{Convergence of the discrete model towards the continuous one}\label{sec: dtoc}

We first notice that the discrete problem \eqref{eq:elli} is a kind of \emph{non conventional} mass-lumped version of a finite element discretization of the continuous one \eqref{eq : continue}. In order to clarify this statement, we introduce the finite element setting. Let $\mathcal{V}_h \subset \mathcal{V}$ the space of continuous, piecewise linear functions $g$ on the one-dimensional partition $T_h = \{y_1, \, \cdots , \, y_{N+1}\}$ of $(0, \Lambda) $, and that verify the Dirichlet boundary condition $g(\Lambda)=0$. Let $\{\Phi_j\}_{j=1,\,N}$ be the standard basis for $\mathcal{V}_h$ consisting of the hat functions defined by $\Phi_j(y_k) = \delta_{j,k}$ for $1\leq j,k\leq N$. 

Let $\ell_h \in \mathcal{V}_h$ be the continuous, piecewise linear function such that for $1 \leq j \leq N+1$, $t>0$, $\ell_h(y_j,t) = \ell_{j}(t)$. Using the basic semi-discrete Galerkine method would lead to the discretization of \eqref{eq:discrGalerkine} in the matrix form: 
\begin{equation}
   \frac{d (M_h L_h)}{dt} + K_h L_h = \tilde f(t),
   \label{eq:discrGalerkine}
\end{equation}
with $L_h(t) = (\ell_1(t), \, \cdots , \, \ell_{N}(t))^T$. Similarly, $\tilde f = (-\frac{\Lambda}{2} \dot L_0, \, 0, \, \cdots, \, 0)$,  $(M_h)_{i,j} = \int\limits_0^{\Lambda} \Phi_i(y) \Phi_j(y)dy$ and $(K_h)_{i,j} = \kappa (\Phi_i,  \Phi_j)$, where $\kappa$ is defined in equation \eqref{eq: Abilin}.

Computing explicitly the coefficients of the matrices $K_h$ and $M_h$ gives
\begin{equation*}
    (K_h)_{ij} = \left\{
    \begin{array}{ll}
    -\Lambda^2K/h & \mbox{for } |i-j|=1,\\
    2\Lambda^2K/h & \mbox{for } i=j\geq 2\,,\\
    \Lambda^2K/h + \Lambda K\tilde{a}/(2a_1) & \mbox{for } i=j=1\,,
    \end{array}
    \right.
\end{equation*}
and
\begin{equation*}
    (M_h)_{ij} = \left\{
    \begin{array}{ll}
    h/6 & \mbox{for } |i-j|=1,\\
    2h/3 & \mbox{for } i=j\geq 2\,,\\
    h/3 & \mbox{for } i=j=1\,.
    \end{array}
    \right.
\end{equation*}
The key observation is that Eqs. (\ref{eq:elli}) and (\ref{eq: l2}) are nothing but a mass-lumped discretization of \eqref{eq : continue} where the mass matrix $M_h$ has been replaced by the diagonal version
\begin{equation*}
    \widetilde{M}_{h}  = \begin{pmatrix}
h  &  & 0\\
 & \ddots  & \\
0 &  & h
\end{pmatrix}\,.
\end{equation*}

Hence, $\ell_h$ actually solves
\begin{equation}
   \frac{d (\widetilde{M}_h L_h)}{dt} + K_h L_h = \tilde f(t)\,,
   \label{eq:mtilde}
\end{equation}
together with the initial condition 
\begin{equation}
    \ell_h(0) = \ell_{0,h} \in \mathcal{V}_h\,.
    \label{eq:mtildeIC}
\end{equation}
The classical mass-lumped method, on the other hand, would have consisted in replacing the tridiagonal mass matrix $M_h$ by a diagonal matrix $\bar{M}_h$ using an integration formula on the vertices of the partition. Namely, using the trapezoidal formula  $\int\limits_0^{\Lambda} g \sim \big(\frac{1}{2}g(y_1) + \sum\limits_{j=2}^{N} g(y_j) + \frac12 g(y_{N+1})\big)h= \big(\frac{1}{2}g(y_1) + \sum\limits_{j=2}^{N} g(y_j)\big)h$, for a function $g$ satisfying $g(\Lambda)=0$ leads to the mass-lumped matrix 

\begin{equation}
    \bar{M}_{h}  = \begin{pmatrix}
h/2  & & & 0\\
 & h   & &\\
  & & \ddots  & \\
0 & & & h 
\end{pmatrix}
\label{eq: MLThomee}
\end{equation}
which differs from $\widetilde{M}_{h}$.

We shall then study the ODE \eqref{eq:mtilde}, \eqref{eq:mtildeIC} using the method presented in \cite{thomee_galerkin_2006} which provides us with a convergence result for the mass-lumped method with $\bar{M}_h$.\\

We introduce the two following inner products on $\mathcal{V}_h$ associated with
$\bar{M}_h$ and $\widetilde{M}_h$ respectively. Namely, for $(u_h,v_h)\in \mathcal{V}_h$
\[ \langle u_h,v_h \rangle_h = \frac{h}{2}u_h(y_1)v_h(y_1) + h\sum\limits_{j=2}^{N} u_h(y_j)v_h(y_j)\]
and
\[(u_h,v_h)_{h} = h\sum\limits_{j=1}^{N} u_h(y_j)v_h(y_j)\,.\]
We also call $\Vert\cdot\Vert_h$ the norm associated to $(\cdot,\cdot)_h$, while the $L^2$ norm and inner products are denoted by $\Vert\cdot\Vert$ and $(\cdot,\cdot)$ respectively. Gerschg\"orin Theorem applied to $M_h$ shows the equivalence of the norms $\Vert\cdot\Vert$ and $\Vert\cdot\Vert_h$ on $\mathcal{V}_h$ \emph{uniformly} in $h$, and, more precisely, we have the estimate, valid for all $v_h\in \mathcal{V}_h$
\[\frac16(v_h,v_h)_h \leq (v_h,v_h) \leq (v_h,v_h)_h\,,\]
from which we also deduce
\begin{equation}
    hv_h(y_1)^2 \leq \Vert v_h\Vert_h^2 \leq 6\Vert v_h\Vert^2\,.
    \label{eq:v_h(y_1)}
\end{equation}

Finally, we introduce, for $u_h,v_h\in \mathcal{V}_h$, $\delta_h (u_h,v_h) = (u_h,v_h)_{h} - (u_h,v_h)$ the quadrature error.

\begin{lemma}
Let $u_h,v_h \in \mathcal{V}_h$. We have, for $h$ sufficiently small:
\begin{align}
    \vert \delta_h (u_h,v_h) \vert &\leq C h \Vert \partial_y u_h \Vert \Vert \partial_y v_h \Vert,
    \label{eq:44}\vspace{6pt}\\
    \vert \delta_h (u_h,v_h) \vert &\leq C \sqrt{h} \Vert \partial_y u_h \Vert \Vert v_h \Vert \label{eq:45}
\end{align}
\label{lemma151}
for a constant $C$ that does not depend on $u_h$, $v_h$ or $h$.
\end{lemma}

\begin{proof}
In all what follows, $C$ denotes a constant that may vary from line to line, being always independent of $h$. Let $u_h,v_h \in \mathcal{V}_h$.
We write $\vert \delta_h(u_h,v_h)\vert \leq \vert(u_h,v_h)_{h}-\langle u_h,v_h \rangle_h\vert +\vert \langle u_h,v_h\rangle_h-(u_h,v_h)\vert$. Thomée \cite{thomee_galerkin_2006} provides us with an estimate of the error between $ \langle u_h,v_h\rangle_{h}$ and $(u_h,v_h)$,  namely, 
\[\vert \langle u_h,v_h \rangle_h - (u_h,v_h) \vert \leq C h^2 \Vert \partial_y u_h\Vert \Vert \partial_y v_h \Vert\] and 
\[\vert \langle u_h,v_h\rangle_h - (u_h,v_h) \vert \leq C h \Vert \partial_y u_h\Vert \Vert v_h \Vert\]
for some constant $C>0$ that does not depend on $u_h$, $v_h$ or $h$. The latter estimate is obtained by an inverse inequality. 

It remains to estimate the term $\tilde \delta_h (u_h,v_h) =  (u_h,v_h)_{h} -  \langle u_h,v_h\rangle_h $.

We notice that:
\begin{align}
    \vert \tilde \delta_h (u_h,v_h) \vert &= \frac{h}{2}\vert u_h(y_1)v_h(y_1)\vert \label{eq:epsh}\vspace{6pt}\\
    &= \frac{h}{2}\left\vert \int_0^\Lambda \partial_y  u_h(y)\,dy\right\vert\,\left\vert \int_0^\Lambda \partial_y v_h(y)\,dy\right\vert \nonumber\\
    & \leq \frac{h\Lambda}{2} \Vert \partial_y u_h \Vert\Vert \partial_y v_h \Vert\,.
    \label{eq:151}
\end{align}

Similarly, \eqref{eq:epsh} together with \eqref{eq:v_h(y_1)} gives:
\begin{equation}
    \vert \tilde \delta_h (u_h,v_h) \vert \leq C \sqrt{h} \Vert \partial_y u_h \Vert\Vert v_h \Vert.
    \label{eq:151p}
\end{equation}

This yields \eqref{eq:44} and \eqref{eq:45}.
\end{proof}

\begin{theorem}
If $\ell$ and $\ell_h$ are solution to~\eqref{eq:FV} and~\eqref{eq:mtilde}, \eqref{eq:mtildeIC} respectively, and $\ell_0\in H^2((0,\Lambda))$,  we have, for all $t \geq 0$,
\begin{eqnarray*}
    \Vert \ell_h(t) - \ell(t) \Vert & \leq C\Vert \ell_{0,h} - \ell_0\Vert + Ch^2(\Vert \partial_{yy}\ell_0 \Vert + \Vert \partial_{yy}\ell(t) \Vert)  \vspace{6pt}\nonumber \\
     & + \,Ch \left( \int\limits_0^t \Vert \partial_{yt}\ell \Vert^2ds \right)^{1/2}. 
\end{eqnarray*}
\label{thm151}
\end{theorem}

\begin{proof}
Let $R_h$ be the Ritz projector from $\mathcal{V}$ on $\mathcal{V}_h$, associated with $\kappa(\cdot,\cdot)$. Namely, for $g\in  \mathcal{V}$, $R_hg$ is defined by

\[\kappa(R_hg,v_h) = \kappa(g,v_h)\]
for all $v_h\in \mathcal{V}_h$.
We write $\ell_h-\ell=\left(\ell_h-R_h \ell\right)+\left(R_h \ell-\ell\right)=\theta_h+\rho$ (Notice that $\theta_h\in \mathcal{V}_h$). Standard estimations show that $\rho(t)$ satisfies $\Vert R_h\ell - \ell\Vert \leq Ch^2\Vert \partial_{yy}\ell\Vert$. In order to estimate $\theta_h$, we write, for all $\chi_h\in \mathcal{V}$

\begin{eqnarray}
(\partial_t \theta_h, \chi_h)_{h}+\kappa(\theta_h, \chi_h) &=&(\partial_t \ell_{h}, \chi_h)_{h}+\kappa(\ell_h, \chi_h) \nonumber \vspace{6pt}\\
&& -(\partial_t R_h \ell, \chi_h)_{h}-\kappa(R_h \ell,\chi_h) \nonumber \vspace{6pt}\\
&=& (f, \chi_h)\nonumber\\
&& -(\partial_t R_h \ell, \chi_h)_{h}-\kappa(\ell,\chi_h) \nonumber \vspace{6pt}\\
&=&(\partial_t \ell, \chi_h)-(\partial_t R_h \ell, \chi_h)_{h} \nonumber \vspace{6pt}\\
&=&-(\partial_t \rho, \chi_h)\nonumber\\
&&\hspace*{1cm}-\delta_h(\partial_t R_h \ell, \chi_h).
    \label{eq:ritz}
\end{eqnarray}

Setting $\chi_h=\theta_h$, we obtain
\begin{equation*}
    \frac{1}{2} \frac{d}{d t}\|\theta_h\|_{h}^2+\kappa(\theta_h, \theta_h)=-\left(\partial_t \rho, \theta_h\right)-\delta_h\left(\partial_t R_h \ell, \theta_h\right) .
\end{equation*}

Here, we have at once, using Cauchy-Schwarz and Poincaré inequalities:
\begin{eqnarray*}
\left|\left(\partial_t \rho, \theta_h\right)\right| &\leq&\left\|\partial_t (\ell-R_h \ell)\right\|\|\theta_h\| \vspace{6pt}\\
&\leq &C h\left\|\partial_{yt} \ell\right\|\|\theta_h\| \vspace{6pt}\\
&\leq& C h\left\|\partial_{yt} \ell\right\|\|\partial_y \theta_h\|.
\end{eqnarray*}

Using the first equation of Lemma \ref{lemma151}, and the fact that $\Vert \partial_y R_h u \Vert \leq C\Vert \partial_y u\Vert$, we also obtain
\begin{eqnarray*}
\left|\delta_{h}\left(\partial_t R_h \ell, \theta_h\right)\right| &\leq& C h\left\|\partial_{yt} R_h \ell\right\|\|\partial_y \theta_h\|\\ &\leq& C h\left\|\partial_{yt}\ell\right\|\|\partial_y \theta_h\|\,,
\end{eqnarray*}
from which we deduce that
\begin{eqnarray*}
\frac{1}{2} \frac{d}{d t}\|\theta_h\|_{h}^2+\kappa(\theta_h, \theta_h) &\leq& C h\left\|\partial_{yt}\ell\right\|\|\partial_y \theta_h\| \\
&\leq& \kappa(\theta_h, \theta_h)+C h^2\left\|\partial_{yt}\ell\right\|^2\,,
\end{eqnarray*}
using the coercivity of $\kappa(\cdot,\cdot)$ on $\mathcal{V}$.
We therefore infer
\begin{equation*}
\|\theta_h(t)\|_{h}^2 \leq\|\theta_h(0)\|_{h}^2+C h^2 \int_0^t\left\|\partial_{yt}\ell\right\|^2\,ds\,.
\end{equation*}
We now recall that $\|\cdot\|_{h}$ and $\|\cdot\|$ are equivalent norms on $\mathcal{V}_h$, uniformly in $h$, and hence
\begin{equation*}
\|\theta_h(t)\| \leq C\|\theta_h(0)\|+C h\left(\int_0^t\left\|\partial_{yt}\ell\right\|^2 ds\right)^{1 / 2} .
\end{equation*}

Here $\|\theta_h(0)\|=\left\|\ell_{0,h}-R_h \ell_0\right\|$ and
\begin{eqnarray*}
\left\|\ell_{0,h}-R_h \ell_0\right\| &\leq& \left\|\ell_{0,h}-\ell_0\right\| + \left\|\ell_{0}-R_h \ell_0\right\|\\
&\leq& \left\|\ell_{0,h}-\ell_0\right\| + C h^2\|\partial_{yy} \ell_0\|,
\end{eqnarray*}
whence $\theta_h(t)$ is bounded as desired.
\end{proof}

\begin{theorem}
If $\ell$ and $\ell_h$ are solution to~\eqref{eq:FV} and~\eqref{eq:mtilde}, \eqref{eq:mtildeIC}  respectively we have, for $t \geq 0$,
\begin{eqnarray*}
    \Vert \partial_y(\ell_h - \ell)(t) \Vert & \leq&  Ch(\Vert \partial_{yy}\ell_0 \Vert + \Vert \partial_{yy}\ell(t) \Vert)  \\
     && \hspace*{-1.5cm}+ C\Vert \partial_y(\ell_{0,h} - \ell_0)\Vert + C\sqrt{h} \left( \int\limits_0^t \Vert \partial_{yt}\ell \Vert^2ds \right)^{1/2}.  
\end{eqnarray*}

\label{thm152}
\end{theorem}

\begin{proof}
We now set $\chi_h = \partial_t \theta_h$ in equation \eqref{eq:ritz} for $\theta_h$ to obtain:
\begin{equation*}
    \Vert \partial_t \theta_h \Vert_{h}^2 + \frac{1}{2}\frac{d}{dt} \kappa(\theta_h, \theta_h) = -(\partial_t \rho, \partial_t\theta_h) - \delta_h(R_h \partial_t\ell, \partial_t\theta_h).
\end{equation*}

Here, as in the proof of Theorem \ref{thm151}, 
\begin{equation*}
    \vert(\partial_t\rho, \partial_t\theta_h)\vert \leq \Vert \partial_t(\ell - R_h \ell)\Vert\Vert \partial_t\theta_h\Vert \leq C \sqrt{h} \Vert \partial_{yt} \ell \Vert \Vert \partial_t \theta_h \Vert.
\end{equation*}

Further, by the second line of Lemma \ref{lemma151}, 
\begin{eqnarray*}
    \vert \delta_h(\partial_t R_h \ell, \partial_t\theta_h)\vert &\leq& C \sqrt{h} \Vert \partial_{yt} R_h \ell\Vert \Vert \partial_t\theta_h\Vert\\
    &\leq& C\sqrt{h} \Vert \partial_{yt}\ell \Vert \Vert \partial_t\theta_h \Vert.
\end{eqnarray*}

Using again the equivalence between the norms $\Vert \cdot \Vert_{h}$ and $\Vert \cdot \Vert$ on $\mathcal{V}_h$, we conclude:
\begin{eqnarray*}
    \Vert \partial_t \theta_h \Vert_{h}^2 + \frac{1}{2}\frac{d}{dt}\kappa(\theta_h, \theta_h) &\leq& C \sqrt{h}\Vert \partial_{yt} \ell \Vert \Vert \partial_t \theta_h \Vert_{h} \\
    &\leq& \Vert \partial_t \theta_h \Vert^2_{h} + Ch\Vert \partial_{yt} \ell \Vert^2\,,
\end{eqnarray*}
so that, after integration, and using the coercivity of $\kappa(\cdot,\cdot)$ on $\mathcal{V}$
\begin{eqnarray*}
    \Vert \partial_y \theta_h(t) \Vert &\leq& C\Vert \partial_y \theta_h(0) \Vert + C \sqrt{h} \left( \int\limits_0^t \Vert \partial_{yt} \ell\Vert^2 ds\right)^{1/2}\\
    & \leq& \Vert \partial_y (\ell_{0,h} -\ell_0)\Vert + Ch\Vert \partial_{yy}\ell_0\Vert \\
    && \hspace{1cm}+ C\sqrt{h}\left(\int\limits_0^t \Vert \partial_{yt} \ell\Vert^2 ds\right)^{1/2}.
\end{eqnarray*}

This, together with the standard estimate for $\partial_y\rho(t)$ completes the proof.
\end{proof}

We proved the convergence of the discrete $N$-spring swimmer to the continuous model we formally derived in the previous section. Note that we obtain only a first-order (resp. half order) convergence in $L^2$ norm  (resp. $H^1$ norm) while the standard estimations for the mass-lumping method leads to a second-order (resp. first order) convergence . This is due to the Fourier-type boundary condition at $0$ which differs from the classical Dirichlet boundary condition used in \cite{thomee_galerkin_2006}.

\section{Mathematical expression of the displacement}\label{sec:dx}
\subsection{Net displacement of the $N$-spring swimmer}
We seek the swimmer's displacement by looking at the displacement of the first of the largest sphere, meaning we only compute $V_1 = \dot x_1$, and integrate over a period $(0,2\pi/\omega)$. 

Taking into account the hydrodynamic interactions due to the $i^{\text{th}}$-sphere with $i\in \{2,\cdots,N+2\}$ on the first sphere, we have
\[\displaystyle V_1 = \frac{1}{6\pi\mu a_1}f_1^F + \frac{1}{4\pi \mu L_0}f_2^F + \frac{1}{4\pi \mu}\sum_{i=3}^{N+2}\frac{f_i^F}{L_0+L_1+\cdots +L_{i-2}}\]

Using expressions \eqref{eq:forcesel} and \eqref{eq:ff}, we obtain

\begin{align}
\begin{array}{ll}
    V_1 &=\displaystyle  \frac{1}{2}\dot L_0 - \frac{\tilde a}{2 a_1} K \ell_1 - \frac{3a_1 \dot L_0}{4 L_0}\vspace{6pt}\\    
    & \displaystyle \hspace*{1cm}-\frac{3K \tilde a \ell_1}{4 L_0} + \frac{3\tilde a K}{2}\sum\limits_{j=1}^{N}  \frac{\ell_{j}-\ell_{j+1}}{\sum\limits_{i=0}^{j} L_i}\,,
    \end{array}
\end{align}
where we recall that, by convention, $\ell_{N+1}=0$.

Finally, by integrating over one period, and noticing that both $\ell_2$ and $\dot L_1 / L_1$ have a vanishing time-average, we obtain, for any value of $h=\Lambda/N$, the displacement of the corresponding $N$-spring swimmer:
\begin{equation}
    \Delta_h x_1 = \int\limits_0^{2\pi/\omega}\bigg[ -\frac{3K \tilde a \ell_1}{4 L_0}+ \frac{3\tilde a K}{2}\sum\limits_{j=1}^{N}  \frac{\ell_{j}-\ell_{j+1}}{\sum\limits_{i=0}^{j} L_i} \bigg]\,dt
    \label{eq:dx1}
\end{equation}

\subsection{Net displacement of the limit model}

We may find an expression for the displacement of the limit model as $h$ tends to 0, by passing to the limit in the preceding expression.

Indeed, for $h$ and $y$ given, we define $j_h(y)$ the unique integer such that $j_h(y)h\leq y\leq (j_h(y)+1)h$. Then, defining $\chi_h$ the function
\[\chi_h(y,t) = \frac{1}{L_0(t)+ \dots +L_{j_h(y)+1}(t)},\]
we may write 
\begin{eqnarray*}
&&\hspace*{-0.7cm}\displaystyle\int_0^{2\pi/\omega}\sum\limits_{j=0}^{N-1}  \frac{\ell_h(jh,t)-\ell_h((j+1)h,t)}{\sum\limits_{i=0}^{j+1} L_i}\,dt=\\
&&\displaystyle\hspace*{2.2cm}-\int_0^{2\pi/\omega}\int_0^\Lambda \partial_y \ell_h(y,t) \chi_h(y,t)\,dy\,dt\,,
\end{eqnarray*}
where $\ell_h$ is the piecewise  linear  function defined in the previous section. 

Finally, the displacement $\Delta_h x_1$ of the $N$-spring swimmer during one time period can be rewritten as
\begin{eqnarray*}
&&    \Delta_h x_1 = \int\limits_0^{2\pi/\omega} \Bigg[-\frac{3 K \tilde a \ell_h(0,t)}{4 L_0(t)}\\
&&\hspace*{1cm}-  \frac{3 \tilde a K}{2}\int_0^\Lambda \partial_y \ell_h(y,t) \chi_h(y,t)\,dy\,\Bigg]\,dt\,.
\end{eqnarray*}

Now, using the fact that $j_h(y)h\to y$ when $h\to 0$, together with the $L^2$ and $H^1$ convergence of $\ell_h$ to $\ell$, we obtain that, for any $y$ and $t$, 
\begin{align}
\chi_h(y,t) &\displaystyle= \frac{1}{\displaystyle L_0(t) + (j_h(y)+1)h + \frac{h}{\Lambda} \sum\limits_{i=0}^{j} \ell_h(ih,t)}\nonumber\\
&\underset{h\to 0}{\longrightarrow} \frac{1}{\displaystyle L_0(t) + y + \int\limits_0^y \frac{\ell(t)}{\Lambda}} =: \chi(y,t)\nonumber
\end{align}

Moreover $0\leq \chi_h(y,t)\leq \max_{t} \frac{1}{L_0(t)}=\frac{1}{L(1-\tilde{\varepsilon})}$, shows that $\chi_h$ is uniformly bounded. 

Therefore, using dominated convergence theorem, we deduce that $\chi_h$ converges to $\chi$ in $L^2(0,2\pi/\omega;(0,\Lambda))$ as $h$ tends to 0. 

Using the convergence theorems proven in the preceding section, we may pass to the limit $h\rightarrow 0$ in $\Delta_h x_1$, and obtain the following expression for the displacement during one period for the limit model 

\begin{equation}
\begin{array}{cc}
    \Delta x_1 & = \displaystyle \int\limits_{0}^{2\pi /\omega} \int\limits_{0}^{\Lambda} -\frac{3 K \tilde a}{2}\partial_y \ell(y,t)\bigg(L_0(t) + y + \displaystyle\int\limits_0^y \frac{\ell}{\Lambda}\bigg)^{-1}\,dy\,dt\\
    & - \displaystyle\int\limits_{0}^{2\pi /\omega}\frac{3 K \tilde a\ell(0,t) }{4 L_0} dt\,.
\end{array}
\end{equation}

\section{Numerical experiments}\label{sec:num}

In this section, we numerically study the discrete model's convergence towards the continuous one. Then, we investigate the influence of the two parameters $\omega$ and $\tilde \varepsilon$ on the system and on its displacement, while the rest of the swimmer is determined by the values in table \ref{tab:param}. All simulations are achieved using \textsc{Matlab}. We consider here that the default length $L$ of the active arm is small compared to the rest of the swimmer. The first sphere thus acts like the head of a sperm cell, and the active arm like a link between the head and the flagella, which gives a signal so that the rest of the system oscillates.

\begin{table}[ht]
    \centering
    \begin{tabular}{cc}
    \hline
        $\tilde a$ &  $1 \cdot 10^{-5} \, m$\\
        $a_1$ &  $1 \cdot 10^{-5} \, m$\\
        $\Lambda$ & $4 \cdot 10^{-4} \, m$\\
        $L$ & $3 \cdot 10^{-5} \, m$\\
        $\tilde k$ & $1 \cdot 10^{-8} \, Nm^{-1}$\\
         $\mu$ & $8.9 \cdot 10^{-4} \, Pa \, s$\\
          \hline
          \vspace{6mm}
         
    \end{tabular}
    \caption{Values of the parameters used in the numerical simulations, matching those of \cite{montino_three-sphere_2015}. We have taken for $\mu$ the dynamic viscosity of water at $25^{\circ}C$.}
    \label{tab:param}
\end{table}

\subsection{Convergence of the discrete models to the continuous one}\label{sec:cvnum}

We investigate numerically the convergence estimations obtained in section \ref{sec: dtoc}. We recall that the continuous solution $\ell$ solves the heat equation PDE with the Fourier-type boundary conditions~\eqref{eq : continue}. We consider, in this section, periodic forcing for which explicit solutions are given by (\ref{eq:Cont_solexpl1}, \ref{eq:Cont_solexpl2}).

\subsubsection{Convergence of the $N$-spring discrete model}

We recall that the discrete solution $\ell_h$ is the $P^1$ discrete function based on the $(\ell_i)_i$ solution to the $N$-spring ODE system (\ref{eq:elli},\ref{eq: lN+2},\ref{eq: l2}). This discrete system corresponds to a semi-discretization in space of the continuous model, based on a non conventional mass-lumping method. The solution $(\ell_i)_i$ of the discrete problem in the periodic setting is given in equations~(\ref{eq:Disc_solexpl1},\ref{eq:gamma},\ref{eq:ab},\ref{eq: bd}). 

\begin{figure}[ht]
    \centering
    \includegraphics[width = 0.75 \textwidth]{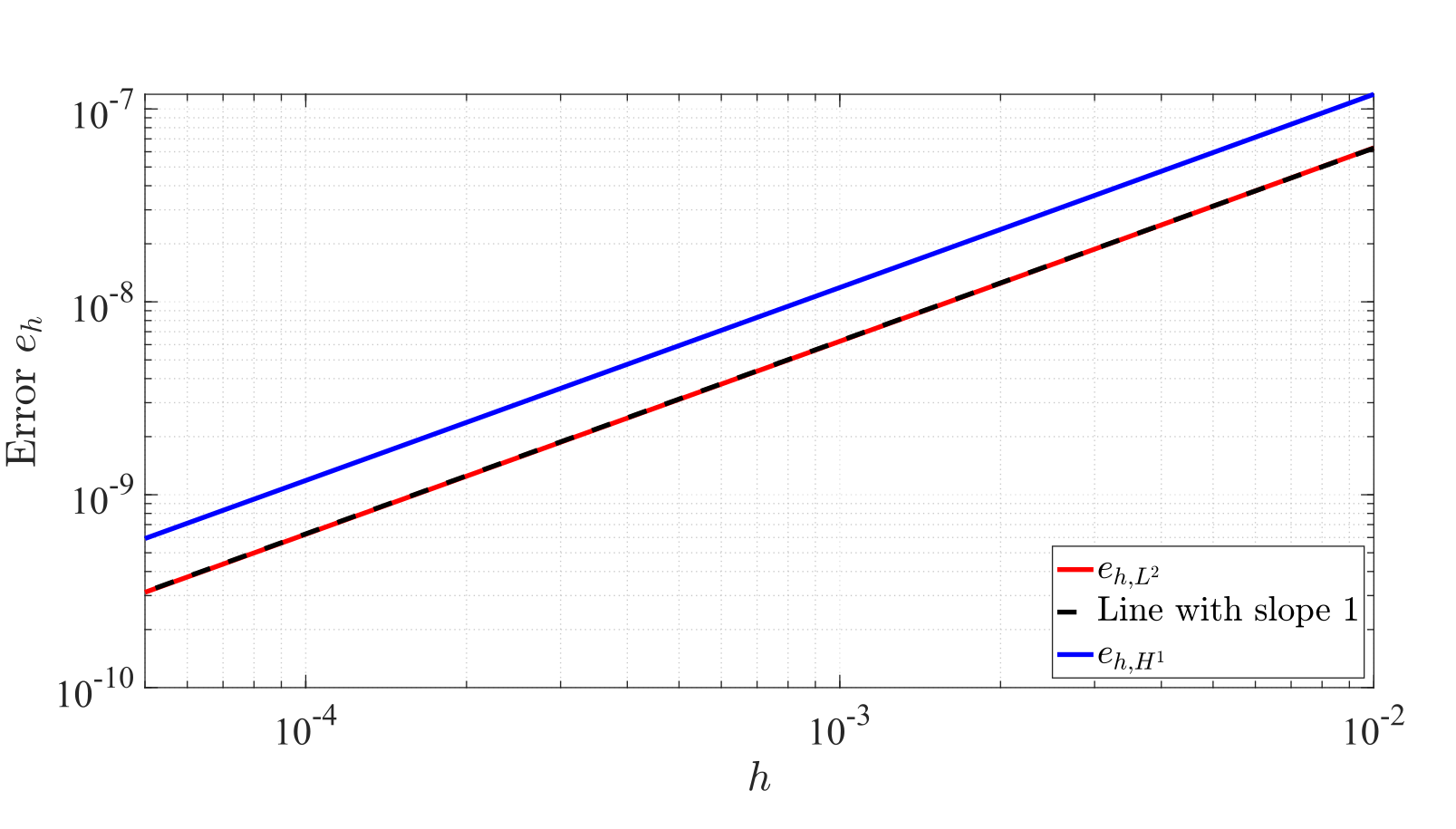}
    \caption{$L^2$ and $H^1$ errors between the $N$-spring discrete model and the continuous one as a function of the number of springs in log scale, in the $(2 \pi/\omega)$-periodic case, for $\tilde \varepsilon = 0.7$ and $\omega = 1 \, rad \cdot s^{-1}$.}
    \label{fig:norm}
\end{figure}

The space step $h$  (or equivalently the number of springs $N$) being given, the discrete error is defined as the error between $\ell_h$ and the $P^1$ interpolation of $\ell$. We plot in figure~\ref{fig:norm}, the $L^2$ (resp. $H^1$) error denoted by $e_{h,L^2}$ (resp. $e_{h,H^1}$). 

We observe that the $L^2$ error converges with order one, as expected from theorem \ref{thm151}. Concerning the $H^1$ error, we observe a superconvergence phenomenon: as the $L^2$ error, it converges at order 1, while theorem~\ref{thm152} predicts a convergence at order 1/2. This can be explained by the regularity of the considered periodic solution.
\\

\subsubsection{Influence of mass-lumping} \label{subsec: MLvscont}

As mentioned earlier, the $N$-spring model turns out to be a discretization in space of the continuous problem~\eqref{eq : continue}, based on an unconventional mass-lumping method. The convergence proof that we proposed in section~\ref{sec: dtoc} is based on the results of Thomée\cite{thomee_galerkin_2006}. He shows that, for a standard mass-lumping discretization, the usual order of convergence for finite elements is obtained: convergence of order 2 for the $L^2$ error and 1 for the $H^1$ error.

We investigate here the influence of the space discretization, by comparing the $N$-spring model~\eqref{eq:mtilde}, solved numerically this time, to the classical mass-lumping method~\eqref{eq: MLThomee} and the standard Galerkin finite element method \eqref{eq:discrGalerkine}.  Again we consider the periodic framework for which the exact solution is available. The time discretization of the three ODE systems is achieved using a Crank-Nicolson scheme for which the time step is chosen to be small enough so that the error due to the time discretization is negligible.

\begin{figure}[ht]
    \centering
    \includegraphics[width = 0.75 \textwidth]{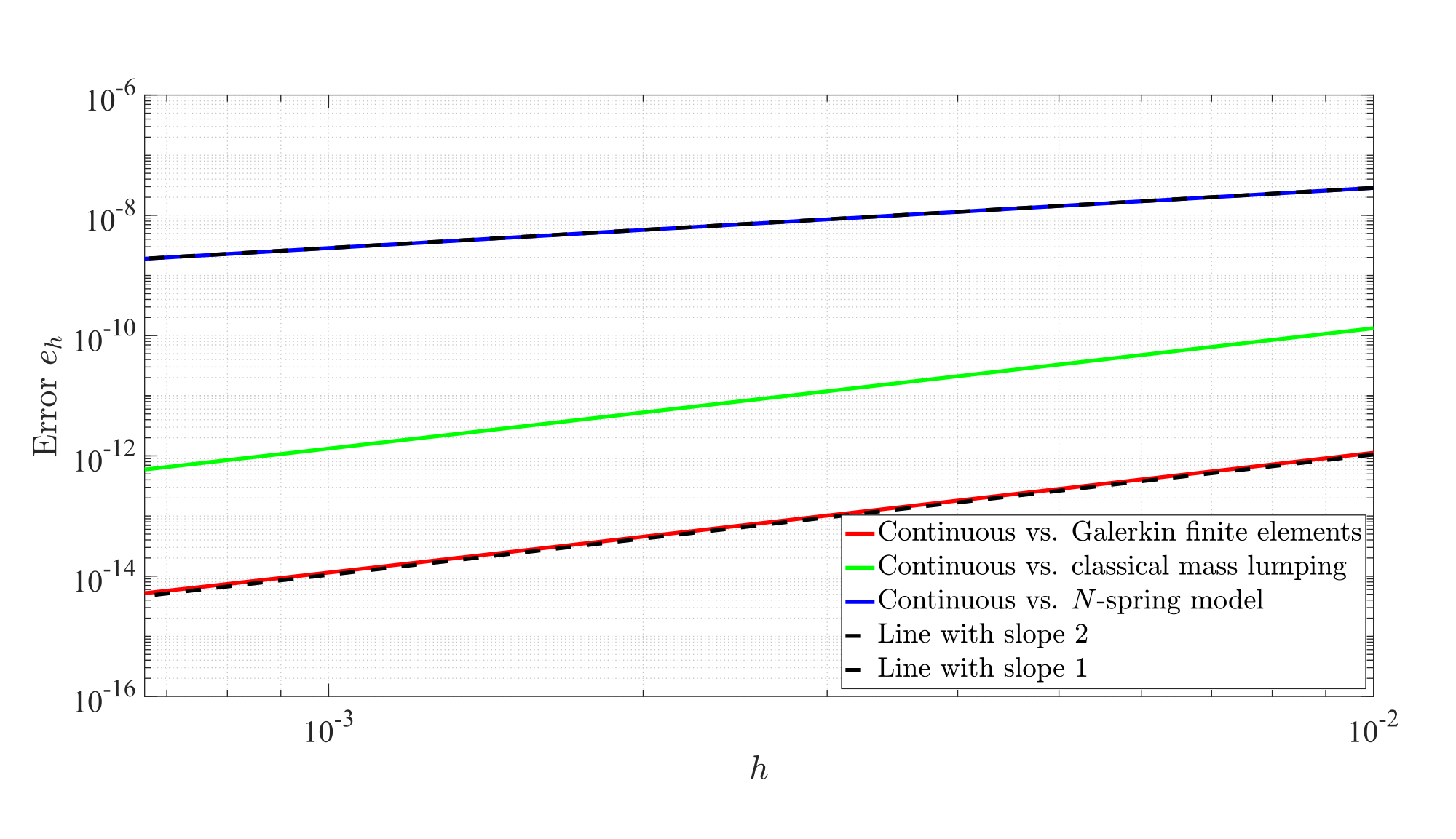}
    \caption{$L^2$ error between the continuous model and our mass-lumping method, as a function of the number of springs, in log scale.}
    \label{fig:L2nonper}
\end{figure}

\begin{figure}[ht]
    \centering
    \includegraphics[width = 0.75 \textwidth]{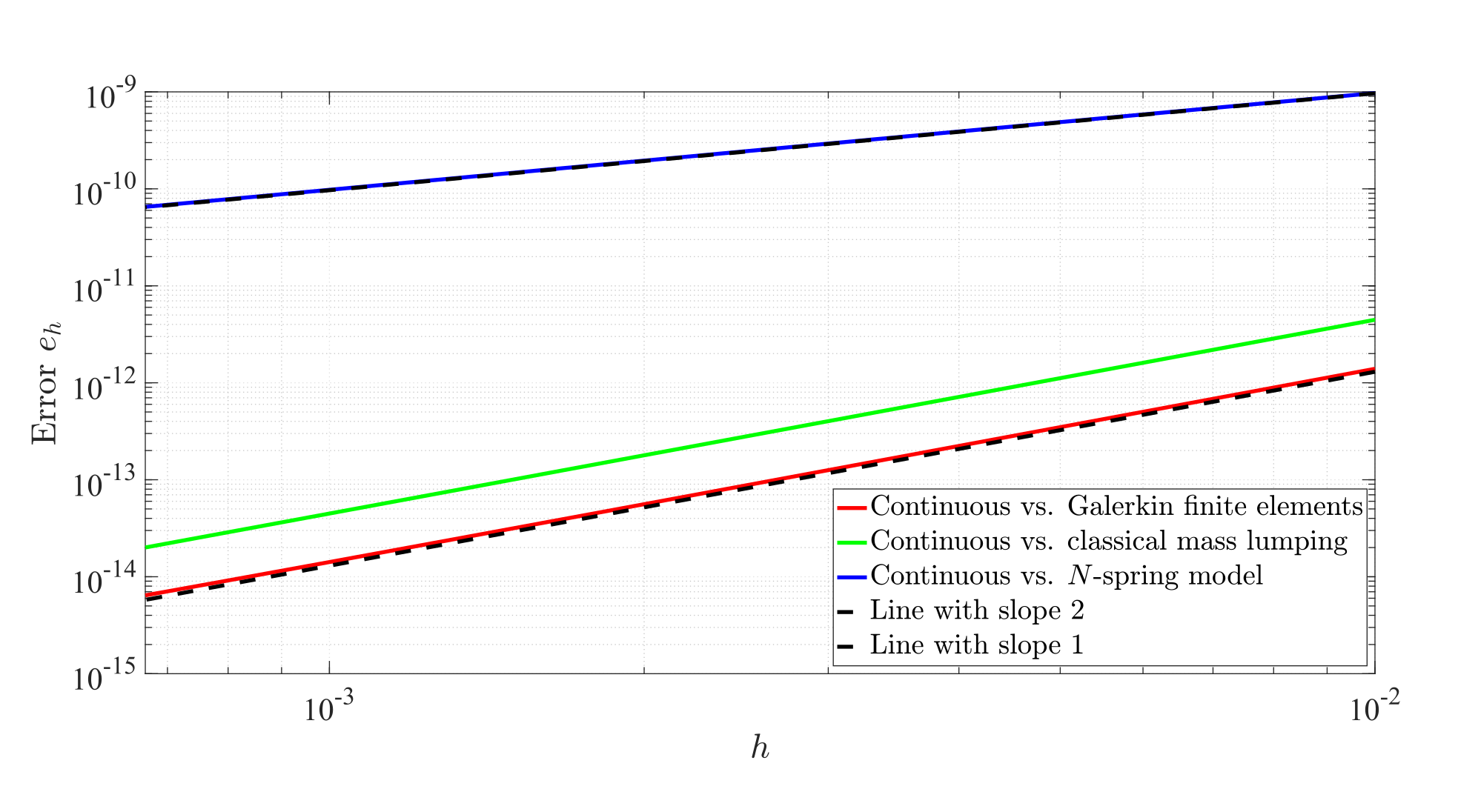}
    \caption{$H^1$ error between the continuous model and our mass-lumping method, as a function of the number of springs, in log scale.}
    \label{fig:H1nonper}
\end{figure}

The corresponding $L^2$ (resp. $H^1$) error is given on figure~\ref{fig:L2nonper} (resp. figure~\ref{fig:H1nonper}). We can see that, as expected, the $L^2$ error converges at order 1 for the $N$-spring model, while it converges at order 2 for both the classical mass-lumping method and the standard Galerking discretization. Again, due to the regularity of the solutions, a super-convergence phenomenon of the $H^1$ error is observed for all three methods: as the $L^2$ error, it converges at order 1 for the $N$-spring model and order 2 for the other two discretizations.

\subsection{Swimming strokes}

In this section, we investigate the swimming ability of the $N$-spring swimmer. The stroke being periodic, we use the explicit solutions given in section~\ref{sec:mvt}. The computations are achieved for $N=2\,000$ springs.

\subsubsection{Deformation of the swimmer}
 
\begin{figure}[ht]
    \centering
\includegraphics[width = 0.75 \textwidth]{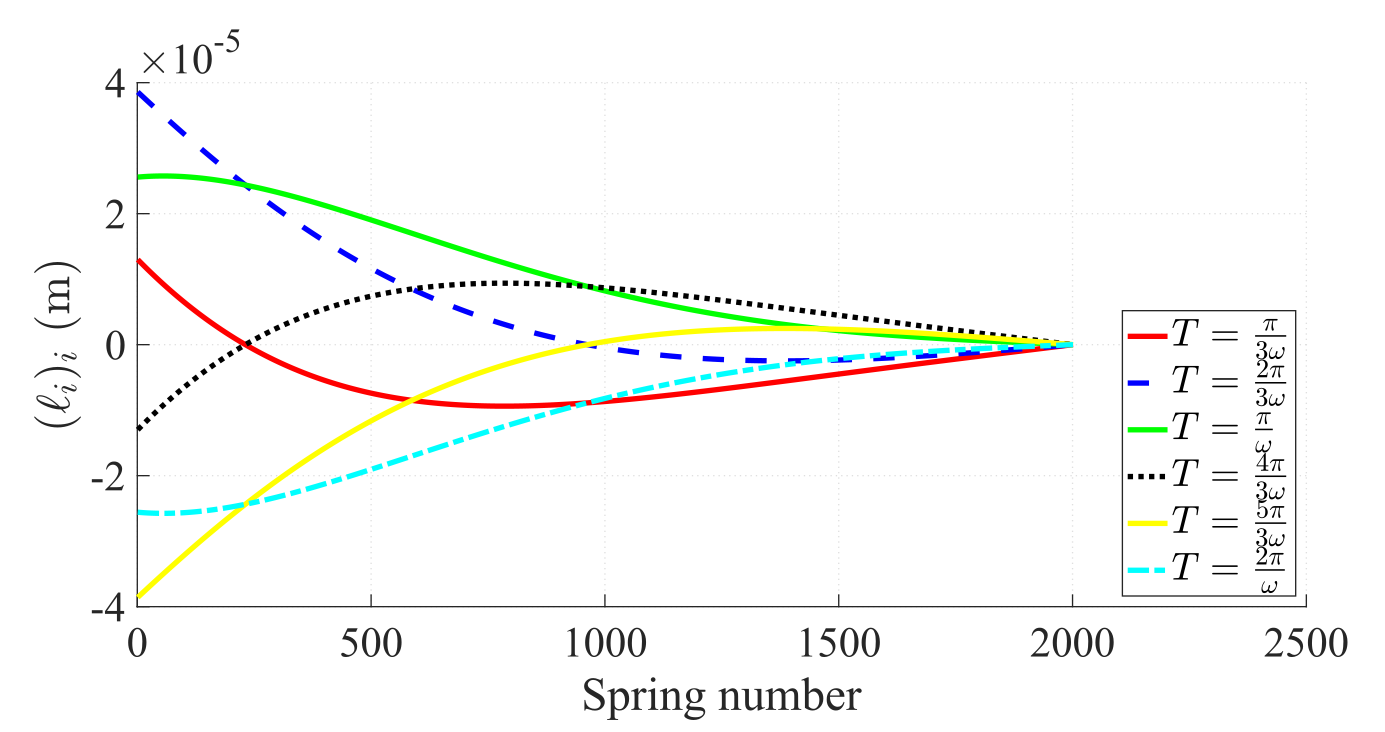}
    \caption{Movement of the whole $2000$-springs swimmer during a full stroke, at different time stamps $T$, for $\omega = 1 \, rad \cdot s^{-1}$ and $\tilde \varepsilon = 0.7$.}
    \label{fig:anim_perl}
\end{figure}

Figure \ref{fig:anim_perl} shows a full stroke of the swimmer, in which we notice that a wave is propagating along its tail. Remember that this wave is a contraction wave along the horizontal tail. This tail appears to be oscillating fairly efficiently for the side close to the head, while the amplitude of the contraction decays considerably on the second half of the tail.

The movement shown corresponds to the stretch of $\ell_j$, and not to the actual deformation which would be $\ell_j/N$, for all $1 \leq j \leq N$. We thus remark that this deformation is relatively small compared to the size of the artificial swimmer, which matches the approximation of small deformations that we made in the first place.

\subsubsection{Displacement}

In this section, we study the influence of the parameters $\tilde \varepsilon$ and $K_{\omega}$ on the swimmer's displacement ~\eqref{eq:dx1}, in order to maximize its absolute value. 

\begin{figure}[ht]
    \centering
\includegraphics[width = 0.75 \textwidth]{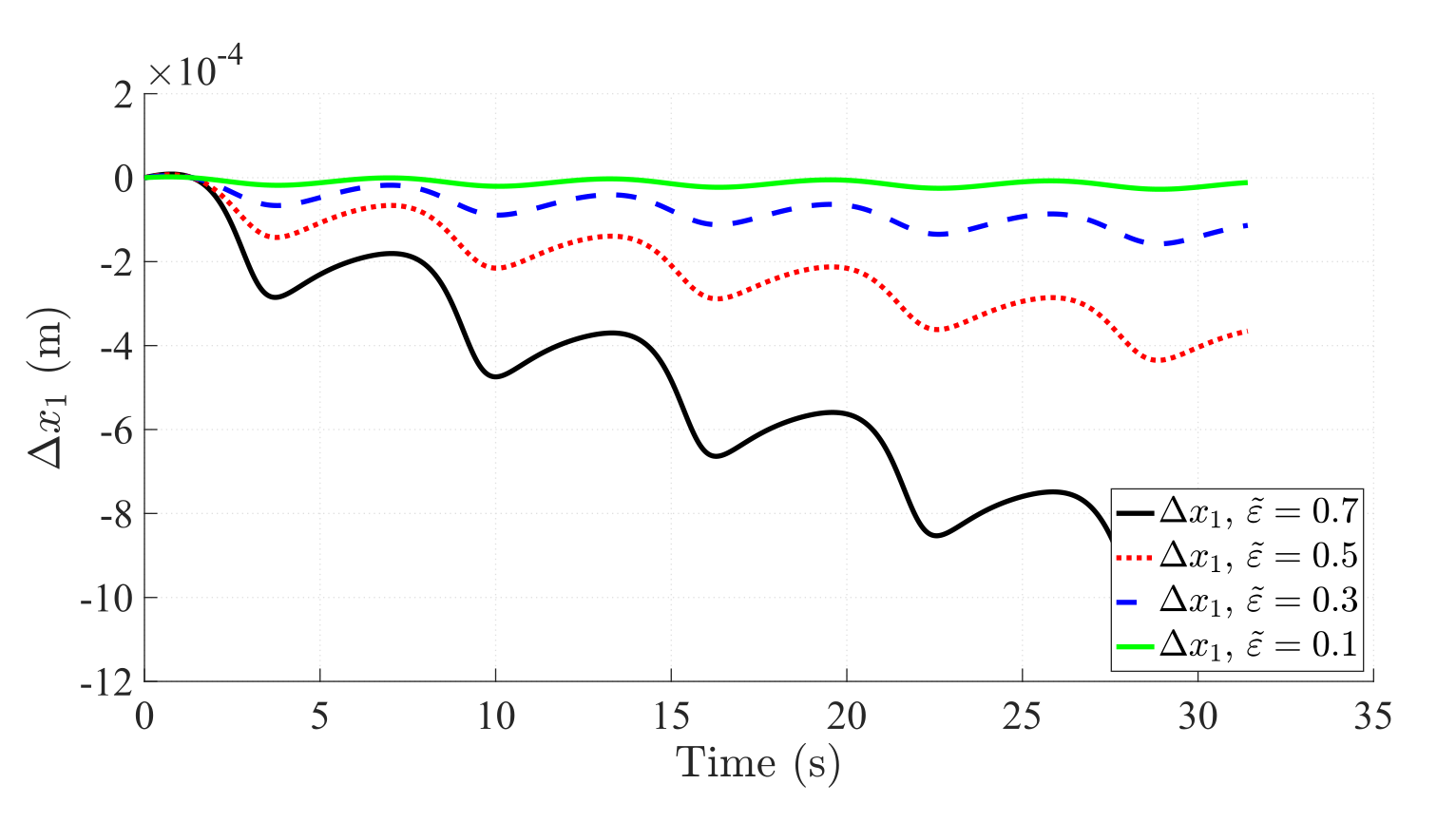}
    \caption{Displacement of the $2000$-spring swimmer against time $t$, for different values of $\tilde \varepsilon$.}
    \label{fig:deltaxt}
\end{figure}

\begin{figure}[ht]
    \centering
\includegraphics[width = 0.75 \textwidth]{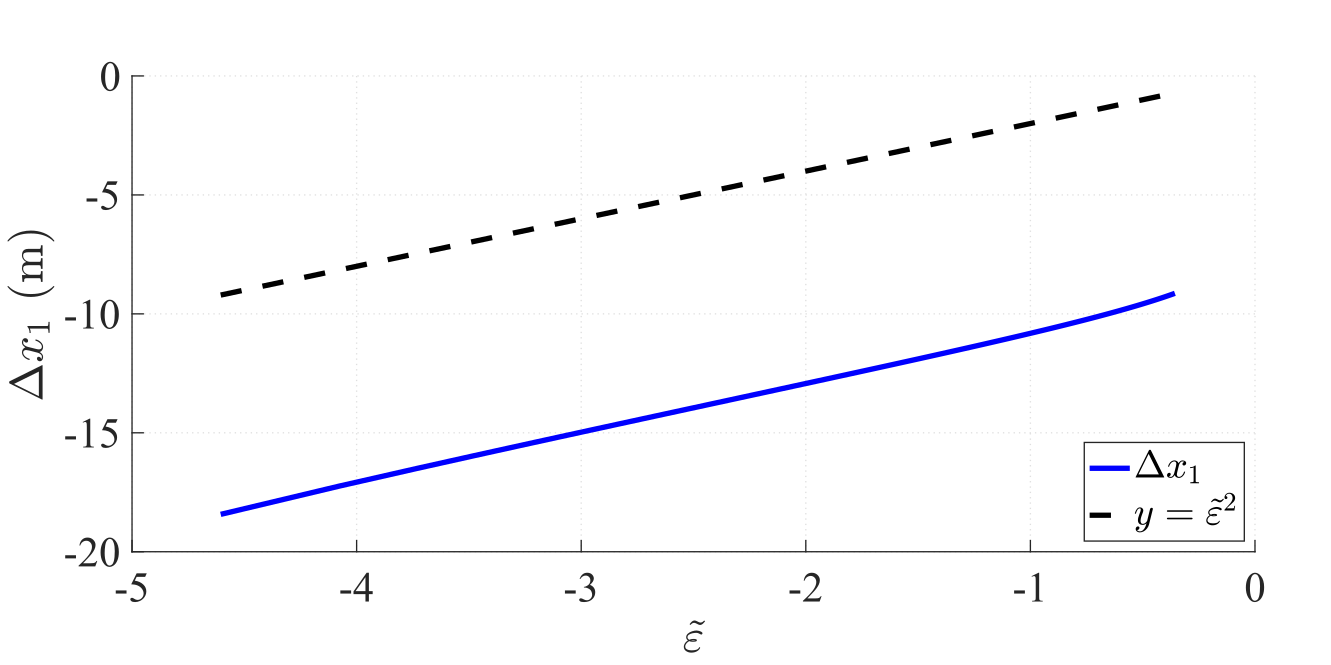}
    \caption{Displacement of the $2000$-spring swimmer depending on $\varepsilon$ for an arm oscillating frequency $\omega = 1 \, rad/s$ compared to $y =\tilde \varepsilon^2 $, in log-scale.}
    \label{fig:deltax}
\end{figure}

In figure~\ref{fig:deltaxt}, we plot the displacement of the swimmer as a function of time, for different values of $\tilde \varepsilon$. The displacement is computed through numerical integration of equation~\eqref{eq:dx1}. The graph shows that the swimmer globally swims backwards, and we recognize the back and forth motion which is characteristic of low Reynolds number artificial swimmers. A larger amplitude $\tilde \varepsilon$ of the forcing leads to a larger displacement and we observe (see figure \ref{fig:deltax}), that $\Delta x_1$ is proportional to $\tilde \varepsilon^2$, which is what is expected in theory (similar behaviors are observed, e.g.,  in \cite{Avron2004,dreyfus_purcells_2005, Avron2005PushmepullyouAE} and explained as the surface of loops in the space of shapes \cite{alouges_optimal_2008}).

 As we want to maximize $\Delta x_1$ while having $\tilde \varepsilon<1$, we choose a fixed value $\tilde \varepsilon = 0.7$ which, although arbitrary, allows for an easier comparison to Montino and DeSimone's results \cite{montino_three-sphere_2015}, as they made a similar parameter choice.

\begin{figure}[!ht]
    \centering
    \includegraphics[width = 0.75 \textwidth]{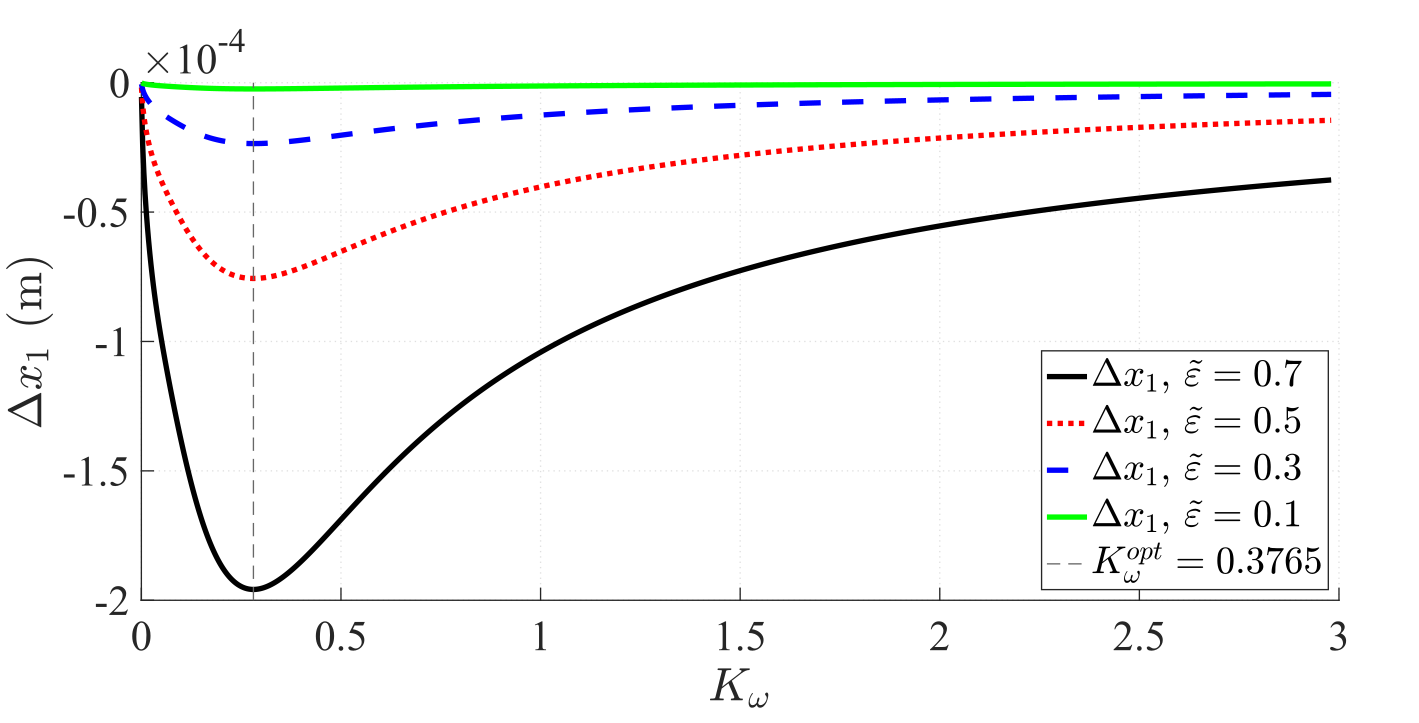}
    \caption{Displacement of the $2000$-spring swimmer depending on $K_{\omega}$, for different values of $\tilde \varepsilon$.}
    \label{fig:dxkom}
\end{figure}

Figure \ref{fig:dxkom} shows $\Delta x_1$ depending on $K_{\omega}$, for different values of $\tilde \varepsilon$. At any fixed $K_{\omega}$, we observe once again that larger $\tilde \varepsilon$ leads to larger $\Delta x_1$. We first observe that, if $K_\omega \to \infty$, the net displacement of the swimmer vanishes. According to the expression of $K_{\omega}$, this is the case for example when $\omega \to 0$: the oscillation disappears, immobilizing the artificial swimmer. This can also happen when $\tilde k \to \infty$: the springs become so rigid that the tail of the swimmer can no longer deform. In that case, the swimmer has only one degree of freedom left to deform and faces Purcell's Scallop theorem's obstruction. Similarly, letting $K_\omega \to 0$ immobilizes the swimmer. An optimal value $K_{\omega}^{\text{opt}}$ for the non-dimensional parameter  is reached between these two limiting cases, in order to maximize the displacement on one time period. According to the figure, $K_{\omega}^{\text{opt}}\simeq 0.3765 $.  A complete mathematical expression of $K_{\omega}^{\text{opt}}$ does not seem available, due to the largely nonlinear nature of the problem contrarily to the final expression obtained in \cite{montino_three-sphere_2015}.
A pair of optimal values for $\omega$ and $k$ to obtain this $K_{\omega}^{\text{opt}}$ are $\omega = 1 \, rads^{-1}$ and $\tilde k \simeq 6.207 \cdot 10^{-8} Nm^{-1}$. Moreover, the expression of $K_{\omega}$ guarantees that $\omega$ must vary proportionally to $\tilde k$ for the pair $(\tilde k,\omega)$ to remain at the optimum.

Indeed, the softer the spring, the slower the first arm needs to oscillate in order to generate a large movement. \\

Looking at the other parameters separately, we can also clearly see from equation \eqref{eq:dx1},  that the displacement depends linearly on $\tilde a$, which is predictable. However, this parameter has a direct consequence on the size of the artificial swimmer and must stay in a reasonable range (in our case no more than $1e-5 \, m$) so that the swimmer stays at microscopic scale.

Finally, we notice that the value of $\Lambda$ and the ratio $a_1/\tilde a$ has little to no influence on the previous analysis. We therefore keep for those parameters values that seem coherent with the scale we are working at, and that match with numerical experiments provided in \cite{montino_three-sphere_2015}.

\section{Conclusion}
We analyzed the dynamics of two low-Reynolds-number swimmers. The first one, which is an extension of \cite{montino_three-sphere_2015}, is made of $N$ passive springs, and the second one is the corresponding limit model with an elastic tail. Both are activated by an active arm that elongates and retracts periodically with amplitude $\varepsilon$ and angular frequency $\omega$.

Noting that the $N$-spring swimmer is a non-conventional mass lumping discretization of the limit model, we proved its convergence, when $N$ tends to infinity, to the continuous model, by extending the results of Thomée~\cite{thomee_galerkin_2006} to the case of a Fourier-type boundary condition.

For both swimmers, a phase difference between the oscillations of the active arm and the tail is created by the interaction between elastic and hydrodynamic forces. Then, both models undergo non-reciprocal shape changes and thus circumvent Scallop Theorem's obstruction \cite{purcell_life_1977}. Numerical simulations indeed show a wave propagating along the swimmers' tails. Similarly to what was shown in \cite{montino_three-sphere_2015}, our models are able to swim but there is no control over the swimming direction.

Then, we focused on computing the net displacement over a time period of the swimmer in both cases, in view of its optimization. We obtain explicit formulae for this displacement as a function of the local elongation during the stroke. We numerically recover the classical back and forth swimming and the second-order scaling of the displacement as a function of the maximum elongation of the forcing active arm. Moreover, we highlight a dimensionless parameter $K_\omega$, driving the movement of the swimmer when its geometry ($\Lambda$, $a$, $a_1$) is given. Some optimal values for this parameter can be estimated by numerical experiments.

Lastly, we noticed that, similarly to the behavior of Machin's swimming rod \cite{machin_wave_1958}, the deformations of  both our swimmers is rapidly attenuating along their passive parts, which suggests that some form of activation is needed in order to mimic the type of behavior observed in nature.

\end{document}